\documentclass[11pt]{article}
\usepackage{lmodern}
\usepackage[T1]{fontenc}
\usepackage{CJKutf8}

\usepackage{graphicx}
\usepackage{float}
\usepackage{amsmath,amsthm,amsfonts,amssymb}
\usepackage{ parskip }
\usepackage{ authblk }
\usepackage{multicol}
\usepackage{caption}
\usepackage{subcaption}
\usepackage[colorlinks=true,linkcolor=blue,citecolor=blue,urlcolor=blue,bookmarks=true]{hyperref}
\usepackage[margin=1.1cm, top=1.5cm, bottom=1.5cm]{geometry}
\usepackage[utf8 ]{inputenc}
\numberwithin{equation}{section}

\begin{document}
	\lineskip 1.5em
	\title{Modified Lotka Volterra  Model with Perspectives of the Piecewise Derivative}
	\author[]{ Atul Kumar }  
	\affil[]{Dayalbagh Educational Institute, Department of Mathematics, India}
	\date{\today}
	\maketitle
	\begin{abstract}
		
		This research presents a concept of piecewise patterns for the different piecewise derivatives using the Lotka Volterra Predator-Prey model. We created numerical solutions using the piecewise derivatives, and known as the Adams-Bashforth approach. The computer results depict piecewise patterns in the real-world behaviors of the Lotka Volterra Predator-Prey model.
		The Lotka-Volterra model investigates the connections between two competing species' abundance and competition. The competitive mechanism is presented and examined, but changes in the abundance of one species are modeled as a function of the abundance of its competitors.
		This idea prompted certain researchers to refer to these mathematical expression as "phenomenological" and to suggest an alternative theoretical structure that focuses particular attention on resources.
		A nonlinear mathematical expression, the Lotka–Volterra model, often known as the Lotka–Volterra predator–prey model, is commonly utilised to discuss the dynamical behaviors of biological systems in which two species connect, one as a predator and the other one as prey. The mathematical expression demonstrates how the populations fluctuate over time.	
		\textbf{Keywords:} Lotka Volterra, Nonlinear System,  Piecewise Derivatives, Classical Differential Equation, Fractional Differential Equation, Stochastic Differential Equation.
	\end{abstract}

	\section{Introduction}
	
	Under several assumptions about the environment and the evolution of predator and prey populations, the Lotka-Volterra model is commonly utilised to assess the connection of two species, called as predator and prey. It should be Keeping in mind that the number of prey directly effects the population of predators. In fact, in the case of limited predators, this model unquestionably assumes a direct connection between the amount of prey and their rate of consumption by a predator \cite{Shim}. For the research, a number of researchers focussed the Lotka Volterra predator-prey model broadly. However, the study of the theory of autocatalytic chemical reaction was first suggested in 1910 by Alfred J. Lotka \cite{Lotka AJ}. He extended the model later, in 1920, to investigate competition between two species \cite{Lotka AJ2}. 
	
	In contrast, fractional calculus approach is also incorporated into a number of predator-prey models since fractional operators are effectual at transferring memory and genetic characteristics of different physical systems \cite{Guimfack, Sekikawa, Danane, Baleanu, Peter, Danane2, k shah, Aldwoah, Atangana2}. For two particular models, the fractional predator-prey model and the fractional rabies model, Ahmed E. et al. \cite{Ahmed E} suggested numerical techniques and addressed the stability of equilibrium points. Din Q. examined and thoroughly studied the equilibrium points of the Lotka-Volterra model \cite{Din Q}. By incorporating a prey refuge into the predator-prey model with fractional derivative, Li HL et al. \cite{Li HL} were able to deduce sufficient requirements for equilibrium point existence.
	A predator-prey model with fractional derivatives was suggested by Ghanbari et al.\cite{Ghanbari B}, along with analysis of their connections. It is significant to notice that the researchers took into account the impact of infection on prey's social behavior as well as its presence in predators. A thorough analysis of a fractional predator-prey system with fear effect was suggested by Yousef et al.\cite{Yousef FB}. The fractional predator-prey system was investigated by El-Saka et al.\cite{El-Saka HA}, who shown the utility of fractional derivative orders as bifurcation parameters in the analysis. The traveling wave solution for the diffusive Lotka–Volterra model was suggested by Tang L et al.\cite{Tang L}. Utilising a non-standard finite difference approach, Eskandari et al.\cite{Eskandari Z} discretized a Lotka–Volterra model.
	
	It is significant to point out that some complex models' behavior may not be adequately represented by the theories and tools now in use. When explaining a technique that transitions from fuzzy to stochastic or from stochastic to power-law, one can utilise the term crossover behavior. The notions of piecewise derivatives operators for modeling complex models for their crossover behaviors have been suggested by several researchers in a number of recent research \cite{Atangana4, Atangana3, Atangana6}. The authors also investigated the cogency of their idea in several models and epidemiological dynamical systems. 
	Other researchers were inspired by the suggested different technique to model novel dynamical systems in order to investigate behaviors of piecewise scenerios, like the Covid-19 spread 	\cite{Atangana3}
	and rhythmic heartbeats \cite{Atangana6}. Author was inspired by this idea to incorporate it into the piecewise Lotka-Volterra Predator-Prey model in order to capture the behaviors of chaos and crossover.
	
	The framework of the study is as follows. In Section 2, we first presented the background information and definitions of piecewise differential operators. Section 3 presents the modified piecewise Lotka-Volterra models. Section 4 includes the numerical approaches and computer results for the considred models at several values of $\alpha$, and concludes with the given information.

	\section{Fundamental Definitions}
	\textbf{Definition 1:} Assume that $\delta \in (0,1)$ and that $Q:[0,P]\rightarrow \mathbb{R}$ is a function. Next, the Caputo fractional derivative\cite{Owolabi} will be given by
	\begin{equation*}	
		_a^CD^{\delta}_tQ(t) = 
		\frac{1}{\Gamma{(1-\delta)}}\displaystyle\int_0^t(t-s)^{-\delta}  { Q'(s)} d{s}.
	\end{equation*}
	
	\textbf{Definition 2:} Let $\delta \in (0, 1]$ and let $Q \in H^{1}(x_{1}, x_{2})$, for $x_{2}>x_{1}$, be a function. Next, the Atangana-Baleanu derivative\cite{Owolabi} is defined as
	
	\begin{equation*}	
		_a^{ABC}D^{\delta}_tQ(t) = 
		\frac{AB(\delta)}{1-\delta}\displaystyle\int_0^t  E_{\delta}\left[-\frac{\delta}{1-\delta}(t-s)^{\delta}\right] Q'(s) ds,
	\end{equation*}
	where $AB(\delta)=1-\delta+\frac{\delta}{\Gamma(\delta)}$ and $E_{\delta}(z)=\displaystyle\sum_{n=0}^{\infty}\frac{z^n}{\Gamma(\delta n + 1)}$ denote the normalization function and the Mittag-Leffler function, respectively.
	
	\textbf{Definition 3:} Let $Q \in H^{1}(x_{1}, x_{2})$ for $x_{2}>x_{1}$ and $\delta \in (0, 1]$. The Caputo–Fabrizio  derivative is provided by \cite{Owolabi}.
	\begin{equation*}	
		_a^{CF}D^{\delta}_tQ(t) = 
		\frac{M(\delta)}{1-\delta}\displaystyle\int_0^t \exp\left[-\frac{\delta}{1-\delta}(t-s)^{\delta}\right]  Q'(s) ds
	\end{equation*}
	using $M(\delta)$ as the normalization function, so that M(0)=M(1)=1.
	
	\textbf{Definition 4:} Using both the traditional and Riemann-Liouville integrals, we establish a piecewise integral for a continuous function $Q$, represented by $_0^{PPL}J^{\delta}_tQ(t)$ \cite{Atangana4}.
	
	\begin{equation}
		_0^{PPL}J^{\delta}_tQ(t) =
		\begin{cases}
			\displaystyle\int_0^{t_{1}} Q(s) d{s} & 0\leq t\leq t_{1} \\
			\frac{1}{\Gamma{(\delta)}}\displaystyle\int_{t_{1}}^t (t-s)^{\delta-1}Q(s)d{s}	& t_{1}\leq t\leq P
		\end{cases}
	\end{equation}

	\textbf{Definition 5: }
	Using the tradional and Caputo-Fabrizio integrals \cite{Atangana4}, we define a piecewise integral for a continuous function $Q$, represented by $_0^{PCF}J^{\delta}_tQ(t)$.
	
	\begin{equation}
		_0^{PCF}J^{\delta}_tQ(t) =
		\begin{cases}
			\displaystyle\int_0^{t_{1}} Q(s) d{s} & 0\leq t\leq t_{1} \\
			\frac{1-\delta}{M(\delta)}Q(t)+\frac{\delta}{M(\delta)}\displaystyle\int_{t_{1}}^t Q(s)(t-s)^{\delta-1}    d{s}	& t_{1}\leq t\leq P
		\end{cases}
	\end{equation}
	
	\textbf{Definition 6: }
	Using the traditional and Atangana-Balaneu integrals \cite{Atangana4}, we define a piecewise integral for a continuous function $Q$, represented by $_0^{PAB}J^{\delta}_tQ(t)$.
	
	\begin{equation}
		_0^{PAB}J^{\delta}_tQ(t) =
		\begin{cases}
			\displaystyle\int_0^{t_{1}} Q(s) d{s} & 0\leq t\leq t_{1} \\
			\frac{1-\delta}{AB(\delta)}Q(t)+\frac{\delta}{AB(\delta)}\frac{1}{\Gamma{(\delta)}}\displaystyle\int_{t_{1}}^t Q(s)(t-s)^{\delta-1}    d{s}	& t_{1}\leq t\leq P
		\end{cases}
	\end{equation}
	%where $	_0^{PAB}J^{\alpha}_tY(t)$ denotes  classical integral on $[0,t_{1}]$ and integral with Mittag-leffler kernel \cite{Atangana4}on  $[t_{1}, \mathcal{T}]$.\\ 
	
	\textbf{Definition 7: } The piecewise derivative of a  function $Q$ defined with a power-law kernel is how  is defined\cite{Atangana4}. 
	\begin{equation}\label{p_caputo}
		_0^{PC}D^{\delta}_tQ(t) = 
		\begin{cases}
			Q'(t) & 0\leq t\leq t_{1} \\
			_{t_{1}}^{C}D^{\delta}_tQ(t)
			& t_{1}\leq t\leq P
		\end{cases}
	\end{equation}

	\textbf{Definition 8: } For every $x_{2}>x_{1}$, let $\delta \in (0, 1]$ and $Q \in H^{1}(x_{1}, x_{2})$. With the help of traditional derivative and the Caputo-Fabrizio derivative, we define a piecewise derivative of $Q$, denoted by $_0^{PCF}D^{\delta}_tQ(t)$ \cite{Atangana4}. 
	
	\begin{equation}\label{p_caputo_fabrizio}
		_0^{PCF}D^{\delta}_tQ(t) = 
		\begin{cases}
			Q'(t) & 0\leq t\leq t_{1} \\
			_{t_{1}}^{CF}D^{\delta}_tQ(t)
			& t_{1}\leq t\leq P
		\end{cases}
	\end{equation}
	
	\textbf{Definition 9:} 
	For every $x_{2}>x_{1}$, let $\delta \in (0, 1]$ and $Q \in H^{1}(x_{1}, x_{2})$. With the help of the traditional derivative and the Atangana-Balaneu derivative in the sense of Caputo \cite{Atangana4}, we create a piecewise derivative of $Q$, represented by $_0^{PAB}D^{\delta}_tQ(t)$.
	\begin{equation}\label{p_atangana_balaneu}
		_0^{PAB}D^{\delta}_tQ(t) = 
		\begin{cases}
			Q'(t) & 0\leq t\leq t_{1} \\
			_{t_{1}}^{ABC}D^{\delta}_tQ(t)
			& t_{1}\leq t\leq P
		\end{cases}
	\end{equation}

	\section{ Lotka Volterra model \cite{Agarwal RP}}	
	The Lotka Volterra prey-predator model shown in this section demonstrates the connections of two species, one acting as a prey and the other one as predator, utilising several patterns. We suppose that such model represents the connections of two species that demonstrates patterns for the same, maintaining universality.  Additionally, several patterns including classical mechanical processes, nonlocal processes, randomness, and their permutations shall be investigated. The utility of such piecewise differential operators in the study of real-world issues that exhibit crossover characteristics must be mentioned. To demonstrate crossover behaviors between two species connecting with distinct patterns, we therefore established randomness and three distinct piecewise differential operators into the Lotka Volterra model.
	
	\subsection*{Case 1: Mathematical model using Power-law kernel \cite{Atangana3}} 
	By applying equation \eqref{p_caputo}, a piecewise system is defined as following
	\begin{align}\label{c_power}
		\begin{cases}
			&\frac{dQ_{i}(t)}{dt} = e(t,Q_{i}), 0 \le t\le P_{1} ,  \\
			& Q_{i}(0)=Q_{i,0}, i=1,2,...,k  \\
			& _{P_{1}}^CD^{\delta}_tQ_{i}(t) = e(t,Q_{i}), P_{1} \le t\le P_{2}  \\
			& Q_{i}(P_{1})=Q_{i,1}, i=1,2,3...k  \\
			&dQ_{i}(t)=e(t,Q_{i})dt+\sigma_{i}Q_{i}d \mathbb{W}_{i}(t),  P_{2} \le t\le P  \\
			& Q_{i}(P_{2})=Q_{i,2}, i=1,2,3...k  \\
		\end{cases}
	\end{align}
	The piecewise Lotka Volterra Prey Predator Model with power-law, which is the mathematical expression of \eqref{c_power}, will be represented by
	\begin{align}
		\begin{cases}
			& \frac{d x(t)}{dt} = x(r-\lambda_{1}x-\lambda_{2}y) \\
			& \frac{d y(t)}{dt} =y(-\lambda_{4}+\lambda_{3}x) \\
			& _{P_{1}}^{C}D^{\delta}_t x(t) =x(r-\lambda_{1}x-\lambda_{2}y) \\
			&  _{P_{1}}^{C}D^{\delta}_t y(t)= y(-\lambda_{4}+\lambda_{3}x)  \\  
			& d x(t)= \bigg(x(r-\lambda_{1}x-\lambda_{2}y)\bigg)dt+\sigma_{1}xd \mathbb{B}_{1}(t) \\
			& d y(t)=\bigg(  y(-\lambda_{4}+\lambda_{3}x)\bigg   )dt+\sigma_{2} yd \mathbb{B}_{2}(t)					
		\end{cases}
	\end{align}
	Where $y(t)$ denotes the number of predators, like foxes, and $x(t)$, the number of its prey, like rabbits; and their derivatives, $\frac{dx(t)}{dt}$ and $\frac{d y(t)}{dt}$, represent the growth of their populations against time $t$. in which $r$, $\lambda_{1}$, $\lambda_{2}$, $\lambda_{3}$, and $\lambda_{4}$ are parameters signifying the connection of those two species. Let $\mathbb{B}_{i}(t)$ denotes Brownian Motions where $i=1,2$. Moreover, $\sigma_{i}$ demonstrates $\mathbb{B}_{i}(t)$'s intensity.

	\subsection*{Case 2: Mathematical model using Mittag-Leffler kernel \cite{Atangana3}}
	Using Mittag-Leffler, a piecewise system is defined as follows by using equation \eqref{p_atangana_balaneu}.
	\begin{align}\label{linear_caputo_mittag}
		\begin{cases}
			&\frac{dQ_{i}(t)}{dt} = e(t,Q_{i}), 0 \le t\le P_{1} ,  \\
			& Q_{i}(0)=Q_{i,0}, i=1,2,...,k  \\
			& _{P_{1}}^{ABC}D^{\delta}_tQ_{i}(t) = e(t,Q_{i}), P_{1} \le t\le P_{2}  \\
			& Q_{i}(P_{1})=Q_{i,1}, i=1,2,3...k  \\
			&dQ_{i}(t)=e(t,Q_{i})dt+\sigma_{i}Q_{i}d \mathbb{W}_{i}(t),  P_{2} \le t\le P  \\
			& Q_{i}(P_{2})=Q_{i,2}, i=1,2,3...k  \\
		\end{cases}
	\end{align}
	The piecewise Lotka Volterra Prey Predator Model with Mittag-leffler, which is the mathematical representation of \eqref{linear_caputo_mittag}, will be denoted by
	
	\begin{align}
		\begin{cases}
			& \frac{d x(t)}{dt} = x(r-\lambda_{1}x-\lambda_{2}y) \\
			& \frac{d y(t)}{dt} =y(-\lambda_{4}+\lambda_{3}x) \\
			& _{P_{1}}^{ABC}D^{\delta}_t x(t) =x(r-\lambda_{1}x-\lambda_{2}y) \\
			&  _{P_{1}}^{ABC}D^{\delta}_t y(t)= y(-\lambda_{4}+\lambda_{3}x)  \\  
			& d x(t)= \bigg(x(r-\lambda_{1}x-\lambda_{2}y)\bigg)dt+\sigma_{1}xd \mathbb{B}_{1}(t) \\
			& d y(t)=\bigg(  y(-\lambda_{4}+\lambda_{3}x)\bigg   )dt+\sigma_{2} yd \mathbb{B}_{2}(t)					
		\end{cases}
	\end{align}
	
	\subsection*{Case 3: Mathematical model using decaying exponential kernel \cite{Atangana3}}
	A piecewise system is defined using fading memory and equation \eqref{p_caputo_fabrizio}.
	\begin{align}\label{linear_caputo_fading}
		\begin{cases}
			&\frac{dQ_{i}(t)}{dt} = e(t,Q_{i}), 0 \le t\le P_{1} ,  \\
			& Q_{i}(0)=Q_{i,0}, i=1,2,...,k  \\
			& _{P_{1}}^{CF}D^{\delta}_tQ_{i}(t) = e(t,Q_{i}), P_{1} \le t\le P_{2}  \\
			& Q_{i}(P_{1})=Q_{i,1}, i=1,2,3...k  \\
			&dQ_{i}(t)=e(t,Q_{i})dt+\sigma_{i}Q_{i}d \mathbb{W}_{i}(t),  P_{2} \le t\le P  \\
			& Q_{i}(P_{2})=Q_{i,2}, i=1,2,3...k  \\
		\end{cases}
	\end{align}
	The piecewise Lotka Volterra Prey Predator Model with fading memory, which is the  mathematical expression of \eqref{linear_caputo_fading}, will be represented by
	\begin{align}
		\begin{cases}
			& \frac{d x(t)}{dt} = x(r-\lambda_{1}x-\lambda_{2}y) \\
			& \frac{d y(t)}{dt} =y(-\lambda_{4}+\lambda_{3}x) \\
			& _{P_{1}}^{CF}D^{\delta}_t x(t) =x(r-\lambda_{1}x-\lambda_{2}y) \\
			&  _{P_{1}}^{CF}D^{\delta}_t y(t)= y(-\lambda_{4}+\lambda_{3}x)  \\  
			& d x(t)= \bigg(x(r-\lambda_{1}x-\lambda_{2}y)\bigg)dt+\sigma_{1}xd \mathbb{B}_{1}(t) \\
			& d y(t)=\bigg(  y(-\lambda_{4}+\lambda_{3}x)\bigg   )dt+\sigma_{2} yd \mathbb{B}_{2}(t)					
		\end{cases}
	\end{align}
	
	\subsection{Stability Analysis of the model }
	A dynamical system's activity may be predicted with the help of the graphical analysis. Nonetheless, there is a approach by which it may be investigated if a model is stable. We shall become proficient in non-linear system analysis. A system of expression for two state variables is as follows:
	\begin{align}
		\begin{cases}
			& _a^{C}D^{\delta}_tx(t) = e_1(x,y) = x(r-\lambda_{1}x-\lambda_{2}y) \\
			&  _a^{C}D^{\delta}_ty(t) = e_2(x,y) = y(-\lambda_{4}+\lambda_{3}x) \\					
		\end{cases}
	\end{align}
	
	To assess the equilibirium points,
	we consider $ _a^{C}D^{\alpha}_tx(t)=0$ and $ _a^{C}D^{\alpha}_ty(t)=0$ then $(x^{eq},y^{eq})=(0,0), (\frac{r}{\lambda_{1}},0), (\frac{\lambda_{4}}{\lambda_{3}},\frac{\lambda_{3}r-\lambda_{1}\lambda_{4}}{\lambda_{3}\lambda_{2}})$ are the equilibirium points.
	For $(x^{eq},y^{eq})=(0,0)$, we find that
	\begin{equation}
		A = 
		\begin{pmatrix}
			r & 0 \\
			0 & -\lambda_{4}  \\
		\end{pmatrix}
	\end{equation}
	then its eigen values are $\nu_{1}=r>0$, $\nu_{2}=-d<0$. Hence equilibirium point  $(x^{eq},y^{eq})=(0,0)$ is unstable.
	
	For $(x^{eq},y^{eq})=(\frac{r}{\lambda_{1}},0)$, we find that
	\begin{equation}
		A = 
		\begin{pmatrix}
			-r & -\frac{r\lambda_{2}}{\lambda_{1}} \\
			0 & \frac{r\lambda_{3}}{\lambda_{1}}-\lambda_{4}  \\
		\end{pmatrix}
	\end{equation}
	then its eigen values are $\nu_{1}=-r<0$, $\nu_{2}=\frac{r\lambda_{3}}{\lambda_{1}}-\lambda_{4} $ if $r\lambda_{3}<\lambda_{1}\lambda_{4}$. Hence equilibirium point  $(x^{eq},y^{eq})=(\frac{r}{\lambda_{1}},0)$ is locally asymptotically stable if $r\lambda_{3}<\lambda_{1}\lambda_{4}$.
	For $(x^{eq},y^{eq})=(\frac{\lambda_{4}}{\lambda_{3}},\frac{r\lambda_{3}-\lambda_{1}\lambda_{4}}{\lambda_{3}\lambda_{2}})$, we find that
	\begin{equation}
		A = 
		\begin{pmatrix}
			-\frac{\lambda_{1}\lambda_{4}}{\lambda_{3}} & -\frac{\lambda_{2}\lambda_{4}}{\lambda_{3}} \\
			\frac{r\lambda_{3}-\lambda_{1}\lambda_{4}}{\lambda_{2}} & 0\\
		\end{pmatrix}
	\end{equation}
	$$\nu_{1}=\frac{-\lambda_{1}\lambda_{4}+\sqrt{(\lambda_{1}\lambda_{4})^2-4\lambda_{3}\lambda_{4}(r\lambda_{3}-\lambda_{1}\lambda_{4})}}{2}$$
	$$\nu_{2}=\frac{-\lambda_{1}\lambda_{4}-\sqrt{(\lambda_{1}\lambda_{4})^2-4\lambda_{3}\lambda_{4}(r\lambda_{3}-\lambda_{1}\lambda_{4})}}{2}$$
	
	\begin{table}[ht]
		\caption{Stabibily analysis} 
		\centering 
		\begin{tabular}{c c c c} 
			\hline\hline 
			Eigen Values & Nature of Points & Stability \\ [0.5ex]
			\hline 
			$\nu_{1}, \nu_{2}>0$ and same signs& node & Unstable(if both $\nu>0$) and Stable(if both $\nu<0$) \\ 
			$\nu_{1}<0$ or $\nu_{2}>0$ and opposite signs  & Saddle point & Unstable  \\
			$\nu_{1}, \nu_{2} \in R$ and $\nu_{1}= \nu_{2}$ & node & Unstable(if both $\nu>0$) and Stable(if both $\nu<0$) \\  
			$\nu_{1}, \nu_{2}=a+ib$; a is non zero & spiral &  Unstable(if Re($\nu)>0$) and Stable(if Re($\nu)<0$) \\
			$\nu_{1}, \nu_{2}=ib$ & centre & Stable but not assymptotic stable  \\ [1ex] 
			\hline 
		\end{tabular}
		\label{table:nonlin} 
	\end{table}
	
	\subsection{Theorem \cite{Agarwal RP}} 
	There exists a constant k > 0 such that
	$|e(t,u)-e(t, v)| \le k|u-v|$, $\forall t \in J,$ and all $u,v \in R$\\
	if $ \frac{kT^{\delta}{(1+\frac{|b|}{|a+b|}})}{\Gamma(\delta+1)}<1$ then the Boundry Value Property that is $_0^{C}D^{\alpha}_t y(t)=e(t,y(t))$ $\forall t \in J=[0,T ]$, $0 < \alpha \le1$, $ay(0) + by(T ) = C
	$ has a unique solution on [0,T ].

	\subsection{Fractional order Lotka Volterra  prey Predator model  }
	We have applied above theorem on the system to check its unique solutions.
	
	\begin{align}
		\begin{cases}
			& _0^{C}D^{\delta}_t x(t) =x(r-\lambda_{1}x-\lambda_{2}y) \\
			&  _0^{C}D^{\delta}_t y(t)= y(-\lambda_{4}+\lambda_{3}x)  \\					
		\end{cases}
	\end{align}
	Applying Lipschitz condition on $e_{1}(t,x,y)$, we have
	$$e_{1}(t,x_{1},y)=x_{1}(r-\lambda_{1}x_{1}-\lambda_{2}y)$$
	$$e_{1}(t,x_{2},y)=x_{2}(r-\lambda_{1}x_{2}-\lambda_{2}y)$$
	$$|e_{1}(t,x_{1},y)-e_{1}(t,x_{2},y))|=|x_{1}(r-\lambda_{1}x_{1}-\lambda_{2}y)-x_{2}(r-\lambda_{1}x_{2}-\lambda_{2}y)|$$
	$$|e_{1}(t,x_{1},y)-e_{1}(t,x_{2},y))|=|(x_{1}-x_{2})r-(x_{1}+x_{2})(x_{1}-x_{2})\lambda_{1}-(x_{1}-x_{2})\lambda_{2}y|$$
	$$|e_{1}(t,x_{1},y)-e_{1}(t,x_{2},y))| \le |(x_{1}-x_{2})||r-(x_{1}+x_{2})\lambda_{1}-\lambda_{2}y|$$
	$$|e_{1}(t,x_{1},y)-e_{1}(t,x_{2},y))| \le |(x_{1}-x_{2})|(r+2\lambda_{1}+\lambda_{2})|$$
	$$|e_{1}(t,x_{1},y)-e_{1}(t,x_{2},y))| \le k|(x_{1}-x_{2})|$$

	where $k = (r+2\lambda_{1}+\lambda_{2})$\\
	If r=0, $\lambda_{1}=0.2$ and $\lambda_{2}=0.1$ then we have k=0.3. We will now verify that the requirement for the suitable value 
	of $\delta \in (0,1]$ with a=1, b=0, $\delta=0.95$ and T=1.
	$$ \frac{kT^{\delta}{(1+\frac{|b|}{|a+b|}})}{\Gamma(\delta+1)}<1$$
	$$ \frac{kT^{\delta}{(1+\frac{|b|}{|a+b|}})}{\Gamma(\delta+1)}=\frac{0.3 \times 1(1+0)}{\Gamma{(0.95+1)}}=\frac{0.3}{\Gamma{(1.95)}}=\frac{0.3}{0.98}=0.306<1$$
	Hence this solution has unique solution.
	
	Similarly, applying Lipschitz condition on $e_{2}(t,x,y) $, we get
	
	$$e_{2}(t,x,y_{1})=y_{1}(-\lambda_{4}+\lambda_{3}x)$$
	$$e_{2}(t,x,y_{2})=y_{2}(-\lambda_{4}+\lambda_{3}x)$$
	$$|e_{2}(t,x,y_{1})-e_{2}(t,x,y_{2})|=|\lambda_{3}xy_{1}-\lambda_{3}xy_{2}-\lambda_{4}(y_{1}-y_{2})|$$
	
	$$|e_{2}(t,x,y_{1})-e_{2}(t,x,y_{2})| \le |(y_{1}-y_{2})||\lambda_{3}x-\lambda_{4}|$$
	$$|e_{2}(t,x,y_{1})-e_{2}(t,x,y_{2})| \le (\lambda_{3}+\lambda_{4})|(y_{1}-y_{2})|$$
	$$|e_{2}(t,x,y_{1})-e_{2}(t,x,y_{2})| \le k|(y_{1}-y_{2})|$$
	
	where $k = (\lambda_{3}+\lambda_{4})$\\
	If $\lambda_{3}=0.1$ and $\lambda_{4}=0.11$ then we have k=0.21. We will now verify that the requirement for the suitable value 
	of $\delta \in (0,1]$ with a=1, b=0, $\delta=0.95$ and T=1.
	$$ \frac{kT^{\delta}{(1+\frac{|b|}{|a+b|}})}{\Gamma(\delta+1)}<1$$
	$$ \frac{kT^{\delta}{(1+\frac{|b|}{|a+b|}})}{\Gamma(\delta+1)}=\frac{0.21 \times 1(1+0)}{\Gamma{(0.95+1)}}=\frac{0.21}{\Gamma{(1.95)}}=\frac{0.21}{0.98}=0.214<1$$
	Hence this solution has unique solution.

	\section{Numerical results}
	In the study, we create numerical approximations for the given problems. In every situation, the Newton polynomial interpolation will be used. The interval $[0,P]$ was first divided as follows:
	\begin{align*}
		0\leq t_{0} \leq t_{1} \leq t_{2}...\leq t_{k_{1}} =\\ P_{1}\leq t_{k_{1}+1} \leq t_{k_{1}+2} \leq t_{k_{1}+3}...\leq t_{k_{2}}= \\P_{2}\leq t_{k_{2}+1} \leq t_{k_{2}+2} \leq t_{k_{2}+3}...\leq t_{k_{3}}=P
	\end{align*}
	Then, using the points $t_{j-2},~t_{j-1},~t_j$, $e(t,X)$ will be approximated by Newton polynomial denoted by $Z(t)$, which will be defined as:
	\begin{align*}\label{Newton_method}
		Z(t)=e(t_{j-1},Q_{j-1})+ \frac{e(t_{j-1},Q_{j-1})-e(t_{j-2},Q_{j-2})}{\Delta t}(t-t_{j-2})+\\ \frac{e(t_{j},Q_{j})-2e(t_{j-1},Q_{j-1})+e(t_{j-2},Q_{j-2})}{2(\Delta t)^2}(t-t_{j-1})(t-t_{j-2})  \nonumber
	\end{align*}

	\subsection{Numerical Simulation Results for first Case}	
	The Newton polynomial-based numerical approximation for the model \eqref{c_power} can be given as
	
	\begin{align}
		\begin{cases}
			Q_{j}^{k_{1}} =& Q_{j}(0)+\frac{1}{12}\sum_{j_{1} = 2}^{k_{1}}\bigg[23e(t_{j}, Q(t_{j}))-16e(t_{j - 1}, Q(t_{j - 1}))+5e(t_{j - 2}, Q(t_{j - 2})) \bigg]* \Delta t, 0 \le t\le P_{1} \\
			Q_{j}^{k_{2}} &= Q_{j}(\mathcal{T}_{1})
			+ \frac{h^\delta}{\Gamma{(\delta)}} \sum_{j_{2} = k_{1}+3}^{k_{2}}\bigg(\frac{e(t_{j - 1}, Q(t_{j -1}))}{1}\bigg[\frac{(k-j)^{\delta}-(k-j-1)^{\delta} }{\delta} \bigg]\\
			&+\frac{e(t_{j - 1}, Q(t_{j - 1}))-e(t_{j - 2}, Q(t_{j - 2}))}{1} \bigg[(k-j+2)\frac{(k-j)^{\delta}-(k-j-1)^{\delta}}{\delta} 
			+\frac{(k-j+1)^{\delta+1}-(k-j)^{\delta+1} }{\delta+1}\bigg]\\
			&+\frac{e(t_{j}, Q(t_{j}))-2e(t_{j-1 }, Q(t_{j-1 }))+e(t_{j_2 }, Q(t_{j_2 }))}{2} \bigg[(k-j+2)(k-j+1)\frac{(k-j+1)^{\delta}-(k-j)^{\delta} }{\delta}\\
			&-(2k-2j+3)\frac{(k-j+1)^{\delta+1}-(k-j)^{\delta+1} }{\delta+1}
			+\frac{(k-j+1)^{\delta+2}-(k-j)^{\delta+2} }{\delta+2}\bigg]\bigg)
			~P_{1} \le t\le P_{2} \\
			Q_{j}^{k_{3}} &= Q_{j}(\mathcal{T}_{2})+ \frac{1}{12} \sum_{j_{3} = k_{2}+3}^{k_{3}}\bigg[23e(t_{j}, Q(t_{j}))- 16e(t_{j - 1}, Q(t_{j - 1})) + 5e(t_{j - 2}, Q(t_{j - 2})) \bigg]*\Delta t+ \\
			&~\sigma_{i} \sum_{j_{3} = k_{2}+3}^{k_{3}}Q(\mathbb{B}_{i}(t)-\mathbb{B}_{i-1}(t)), P_{2} \le t\le P 
		\end{cases}
	\end{align} 
	\subsection*{Numerical simulation of piecewise Lotka Volterra Predator Prey model  for first case:}  
\begin{figure}[H]
		\minipage{0.49\textwidth}
		\centering
		\includegraphics[width=\linewidth]{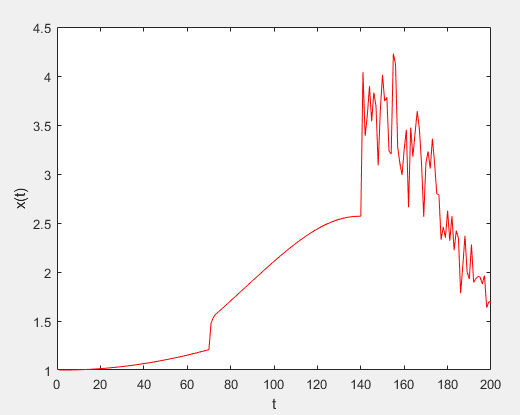}  
		\caption{ $\delta=0.91$.}
		\label{fig:fig1}
		\endminipage\hfill
		\minipage{0.49\textwidth}
		\centering
		\includegraphics[width=\linewidth]{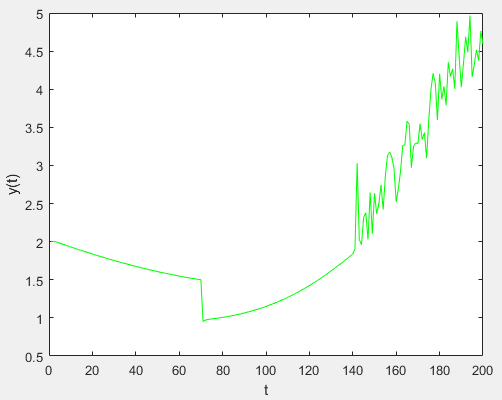}  
		\caption{ $\delta=0.91$.}
		\label{fig:fig2}
		\endminipage\hfill
		\begin{center}
			\minipage{0.63\textwidth}%
			\centering
			\includegraphics[width=\linewidth]{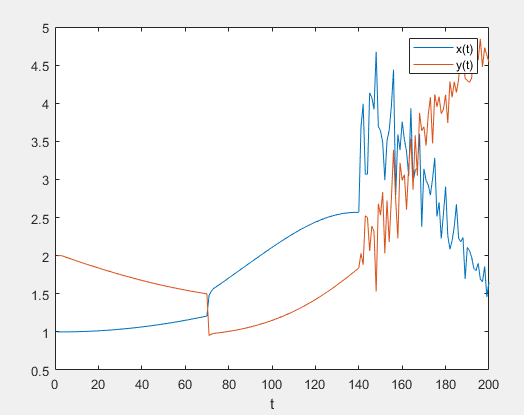}  
			\caption{ $\delta=0.91$.}
			\label{fig:fig3}
			\endminipage
		\end{center}
	\end{figure}

	Through the first case, it presents the numerical solutions to the piecewise Lotka Volterra Predator Prey model. We use the starting values as
	$t(0)=0$,  $h=0.01$, $x(0)=1$, $y(0)=2$, $\lambda_{1}=2$, $\lambda_{2}=1$, $\lambda_{3}=1.5$, $\lambda_{4}=1$, $\sigma_{1}=0.1$ and $\sigma_{2}=0.1$.
	See figures \ref{fig:fig1}, \ref{fig:fig2} and \ref{fig:fig3}.

	\subsection*{Chaotic behavior of numerical simulation of piecewise Lotka Volterra Predator Prey model  for first case: }  
\begin{figure}[H]
		\centering
		\subfloat[for $\delta=0.6$]{\includegraphics[width=0.45\textwidth]{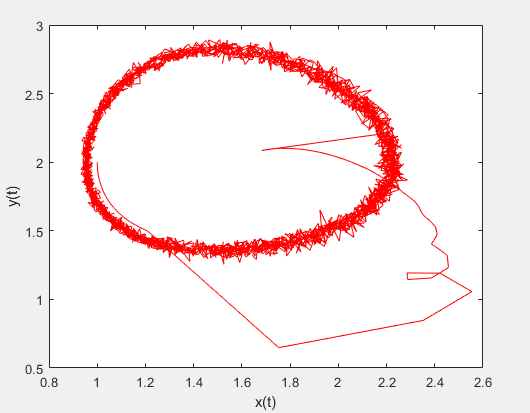}}\hspace{0.5cm}
		%\label{fig:fig4a}
		\subfloat[for $\delta=0.68$]{\includegraphics[width=0.47\textwidth]{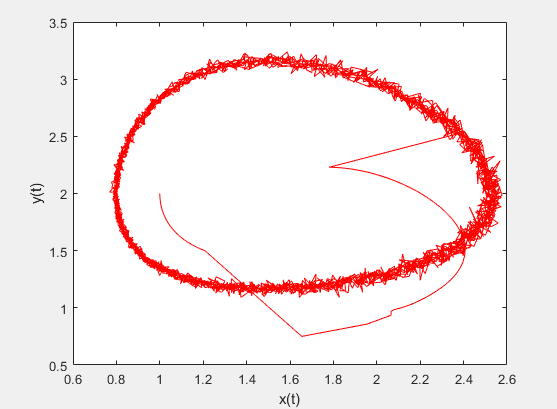}}\\
		%\label{fig:fig4b}
		\caption*{}
		\subfloat[for $\delta=0.85$]{\includegraphics[width=0.46\textwidth]{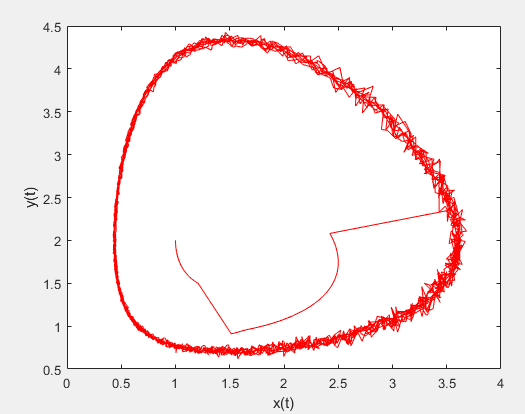}}\hspace{0.5cm}
		%\label{fig:fig4c}
		\subfloat[for $\delta=0.97$]{\includegraphics[width=0.45\textwidth]{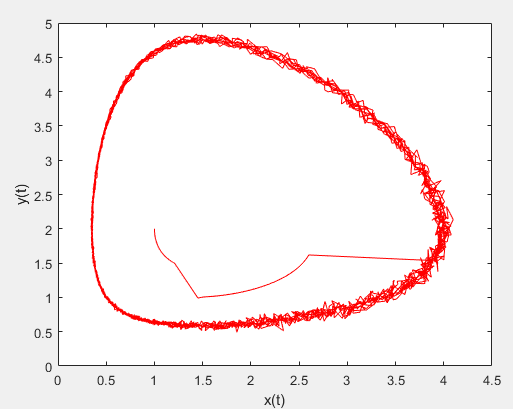}}\\
		%		\caption*{B-picture (a) and (b)}
		\caption{Phase Portrait between $x(t)$ and $y(t)$}
		\label{fig:fig4d}	
	\end{figure}

	Through the first case, it presents the numerical solutions to the piecewise Lotka Volterra Predator Prey model. We use the starting values as
	$t(0)=0$,  $h=0.01$, $x(0)=1$, $y(0)=2$, $\lambda_{1}=2$, $\lambda_{2}=1$, $\lambda_{3}=1.5$, $\lambda_{4}=1$, $\sigma_{1}=0.01$ and $\sigma_{2}=0.02$.
	See figures  \ref{fig:fig4d}(a), \ref{fig:fig4d}(b), \ref{fig:fig4d}(c), and \ref{fig:fig4d}(d).

	\subsection{Numerical Simulation Results for second Case}
	The second case \eqref{linear_caputo_mittag} is applied in this section. The numerical scheme, applying the Newton polynomial, will be provided as
	\begin{align}
		\begin{cases}
			Q_{j}^{k_{1}}&= Q_{j}(0)+\frac{1}{12} \sum_{j_{1} = 2}^{k_{1}} \bigg[23e(t_{j}, Q(t_{j}))- 16e(t_{j - 1}, Q(t_{j - 1})) + 5e(t_{j - 2}, Q(t_{j - 2})) \bigg]*\Delta t, 0 \le t\le P_{1} \\
			Q_{j}^{k_{2}} &= Q_{j}(\mathcal{T}_{1})+
			\frac{1-\delta}{AB(\delta)} g(t_{k}, Q(t_{k}))+\\
			&+\frac{\delta}{AB(\delta) \Gamma(\delta+1)} \frac{h^\delta}{\Gamma{(\delta)}} \sum_{j_{2} = k_{1}+3}^{k_{2}}\bigg(\frac{e(t_{j - 1}, Q(t_{j -1}))}{1}\bigg[\frac{(k-j)^{\delta}-(k-j-1)^{\delta} }{\delta} \bigg]\\
			&+\frac{e(t_{j - 1}, Q(t_{j - 1}))-e(t_{j - 2}, Q(t_{j - 2}))}{1} \bigg[(k-j+2)\frac{(k-j)^{\delta}-(k-j-1)^{\delta}}{\delta} 
			+\frac{(k-j+1)^{\delta+1}-(k-j)^{\delta+1} }{\delta+1}\bigg]\\
			&+\frac{e(t_{j}, Q(t_{j}))-2e(t_{j-1 }, Q(t_{j-1 }))+e(t_{j_2 }, Q(t_{j_2 }))}{2} \bigg[(k-j+2)(k-j+1)\frac{(k-j+1)^{\delta}-(k-j)^{\delta} }{\delta}\\
			&-(2k-2j+3)\frac{(k-j+1)^{\delta+1}-(k-j)^{\delta+1} }{\delta+1}
			+\frac{(k-j+1)^{\delta+2}-(k-j)^{\delta+2} }{\delta+2}\bigg]\bigg),
			P_{1} \le t\le P_{2}\\
			Q_{i}^{k_{3}}&= Q_{j}(\mathcal{T}_{2})+ \frac{1}{12}\sum_{j_{3} = k_{2}+3}^{k_{3}} \bigg[23e(t_{j}, Q(t_{j}))-16e(t_{j - 1}, Q(t_{j - 1})) +5e(t_{j - 2}, Q(t_{j - 2}))\bigg]*\Delta t+ \\
			&~\sigma_{i}\sum_{j_{3} = k_{2}+3}^{k_{3}} Q(\mathbb{B}_{i}(t)-\mathbb{B}_{i-1}(t)), P_{2} \le t\le P
		\end{cases}
	\end{align}
	
	\subsection*{ Numerical simulation of Lotka Volterra Predator Prey model  for  Classical Mittag-Leffler-law-randomness  }  
\begin{figure}[H]
		\minipage{0.49\textwidth}
		\centering
		\includegraphics[width=\linewidth]{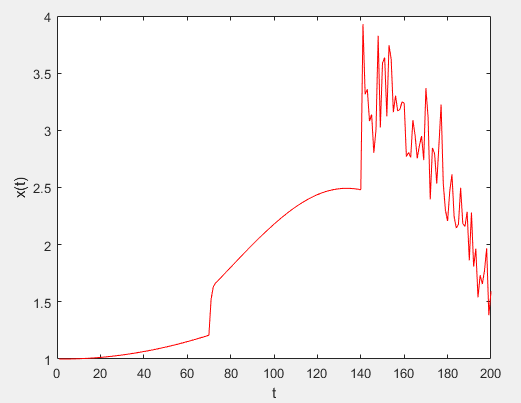}  
		\caption{ $\delta=0.8$.}
		\label{fig:fig5}
		\endminipage\hfill
		\minipage{0.49\textwidth}
		\centering
		\includegraphics[width=\linewidth]{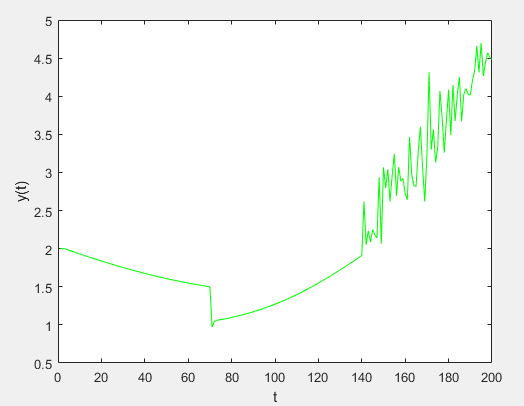}  
		\caption{$\delta=0.8$.}
		\label{fig:fig6}
		\endminipage\hfill
		\begin{center}
			\minipage{0.65\textwidth}%
			\centering
			\includegraphics[width=\linewidth]{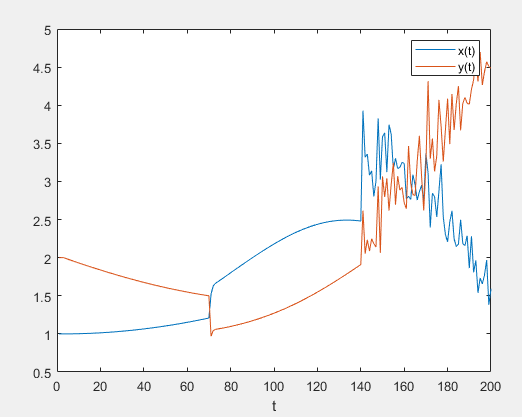}  
			\caption{ $\delta=0.8$.}
			\label{fig:fig7}
			\endminipage
		\end{center}
	\end{figure}

	By applying the Mittag-Leffler-law operator, it presents the numerical solutions to the piecewise Lotka Volterra predator-prey model.The starting circumstances were regarded as
	$t(0)=0$,  $h=0.01$, $x(0)=1$, $y(0)=2$, $\lambda_{1}=2$, $\lambda_{2}=1$, $\lambda_{3}=1.5$, $\lambda_{4}=1$, $\sigma_{1}=0.1$ and $\sigma_{2}=0.11$.
	See figures  \ref{fig:fig5}, \ref{fig:fig6} and \ref{fig:fig7}.

	\subsection*{Chaotic behavior of numerical simulation of piecewise piecewise Lotka Volterra predator prey model for  second case: }  
\begin{figure}[H]
		\centering
		\subfloat[for $\delta=0.75$]{\includegraphics[width=0.45\textwidth]{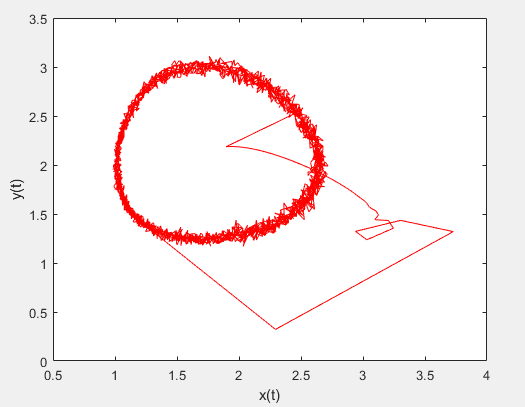}}\hspace{0.5cm}
		%\label{fig:fig8a}
		\subfloat[for $\delta=0.82$]{\includegraphics[width=0.45\textwidth]{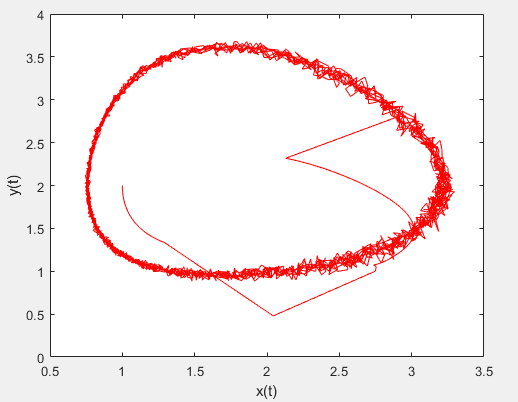}}\\
		\caption*{}
		%\label{fig:fig8b}
		\subfloat[for $\delta=0.89$]{\includegraphics[width=0.47\textwidth]{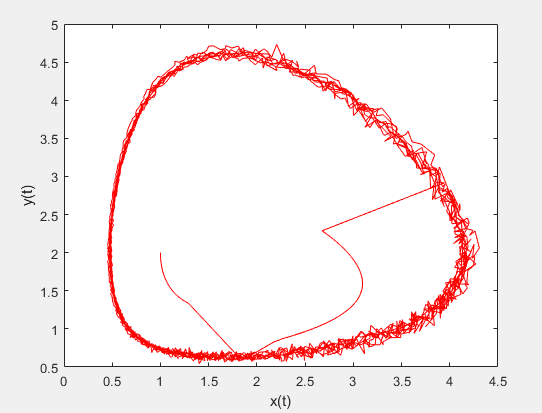}}\hspace{0.5cm}
		%\label{fig:fig8c}
		\subfloat[for $\delta=0.98$]{\includegraphics[width=0.45\textwidth]{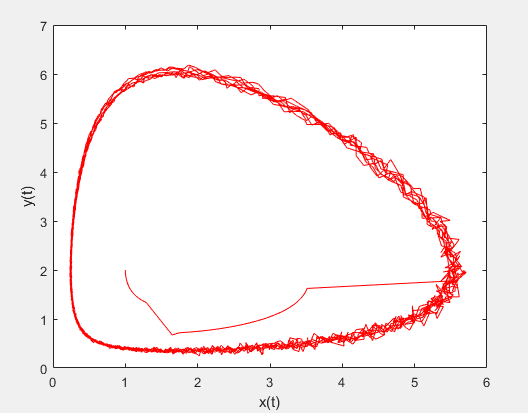}}\\
		%		\caption*{B-picture (a) and (b)}
		\caption{Phase Portrait between $x(t)$ and $y(t)$} 
		\label{fig:fig8d}
	\end{figure}

	By using the Mittag-Leffler-law, it presents the numerical solutions to the piecewise Lotka Volterra predator-prey model.The starting circumstances were regarded as
	$t(0)=0$,  $h=0.01$, $x(0)=1$, $y(0)=2$, $\lambda_{1}=2$, $\lambda_{2}=1$, $\lambda_{3}=1.7$, $\lambda_{4}=1.7$, $\sigma_{1}=0.01$ and $\sigma_{2}=0.02$. 
	See figures  \ref{fig:fig8d}(a), \ref{fig:fig8d}(b), \ref{fig:fig8d}(c), and \ref{fig:fig8d}(d).
	
	\subsection{Numerical Simulation Results for third Case}
	The third case, \eqref{linear_caputo_fading}, is now under investigation. The Newton Polynomial is used to obtain the numerical scheme, which is provided as
	\begin{align}
		\begin{cases}
			Q_{j}^{k_{1}}&= Q_{j}({0})+ \frac{1}{12}\sum_{j_{1} = 2}^{k_{1}} \bigg[23e(t_{j}, Q(t_{j}))- 16e(t_{j - 1}, Q(t_{j - 1})) + 5e(t_{j - 2}, Q(t_{j - 2})) \bigg]*\Delta t, 0 \le t\le P_{1} \\
			Q_{j}^{k_{2}}&= Q_{j}(P_{1})+\frac{1-\delta}{M(\delta)} \bigg[ e(t_{k}, Q(t_{k})) - e(t_{k - 1}, Q(t_{k - 1})) \bigg] + \\ &~\frac{1}{12}\sum_{j_{2} = k_{1}+1}^{k_{2}} \bigg[23e(t_{j}, Q(t_{j})) -16e(t_{j - 1}, Q(t_{j - 1}))+5e(t_{j - 2}, Q(t_{j - 2})) \bigg]*\Delta t, 
			\mathcal{T}_{1} \le t\le P_{2} \\
			Q_{j}^{k_{3}}&= Q_{j}(P_{2}) + \frac{1}{12}\sum_{j_{3} = k_{2}+3}^{k_{3}} \bigg[23e(t_{j}, Q(t_{j}))-16e(t_{j - 1}, Q(t_{j - 1})) + 5e(t_{j - 2}, Q(t_{j - 2})) \bigg]*\Delta t+ \\ &~\sigma_{i}\sum_{j_{3} = k_{2}+3}^{k_{3}}  Q_{i}^{k_{3}}(\mathbb{B}_{i}(t)-\mathbb{B}_{i-1}(t)), \\
			&P_{2} \le t\le P 
		\end{cases}
	\end{align}
	
	\subsection*{ Numerical simulation of  piecewise Lotka Volterra predator prey model for third case}  
\begin{figure}[H]
		\minipage{0.49\textwidth}
		\centering
		\includegraphics[width=\linewidth]{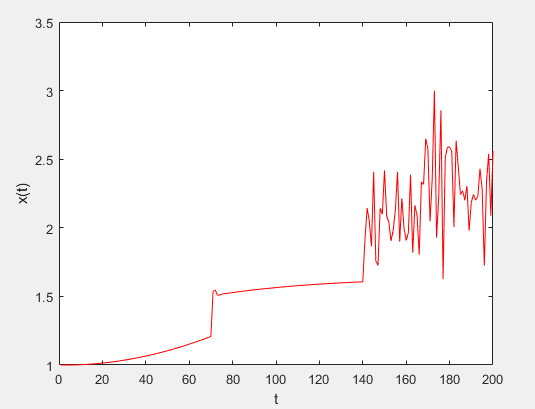}  
		\caption{ $\delta=0.94$.}
		\label{fig:fig9}
		\endminipage\hfill
		\minipage{0.49\textwidth}
		\centering
		\includegraphics[width=\linewidth]{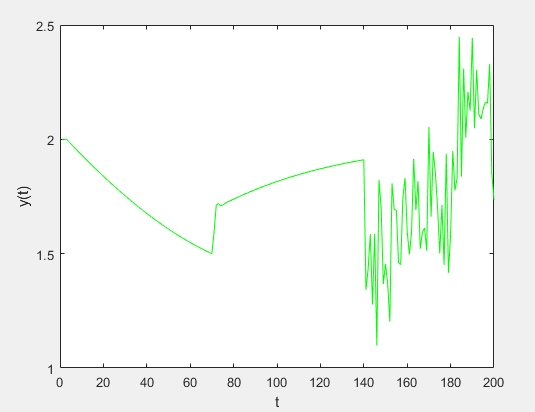}  
		\caption{ $\delta=0.94$.}
		\label{fig:fig10}
		\endminipage\hfill
		\begin{center}
			\minipage{0.7\textwidth}%
			\centering
			\includegraphics[width=\linewidth]{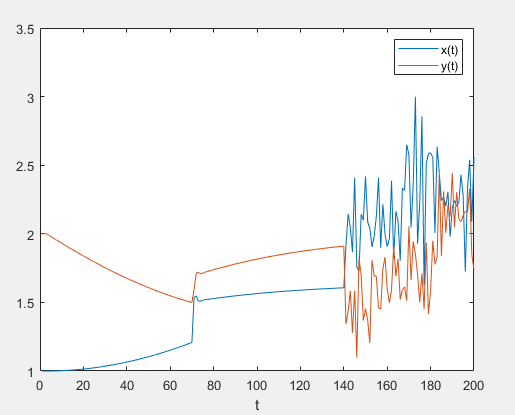}  
			\caption{ $\delta=0.94$.}
			\label{fig:fig11}
			\endminipage
		\end{center}
	\end{figure}

	By using the exponential decay, it presents the numerical solutions to the piecewise Lotka Volterra predator-prey model.The starting circumstances were regarded as
	$t(0)=0$,  $h=0.01$, $x(0)=1$, $y(0)=2$, $\lambda_{1}=2$, $\lambda_{2}=1$, $\lambda_{3}=1.5$, $\lambda_{4}=1$, $\sigma_{1}=0.1$ and $\sigma_{2}=0.11$.
	See figures  \ref{fig:fig9}, \ref{fig:fig10} and \ref{fig:fig11}.

	\subsection*{Chaotic behavior of numerical simulation of  piecewise Lotka Volterra predator prey model for third case}  
\begin{figure}[H]
		\centering
		\subfloat[for $\delta=0.68$]{\includegraphics[width=0.45\textwidth]{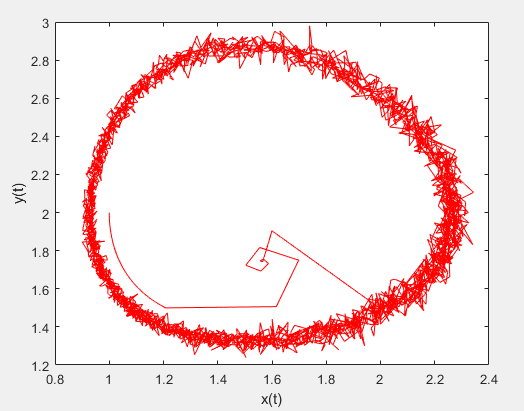}}\hspace{0.5cm}
		\label{fig:fig12a}
		\subfloat[for $\delta=0.78$]{\includegraphics[width=0.45\textwidth]{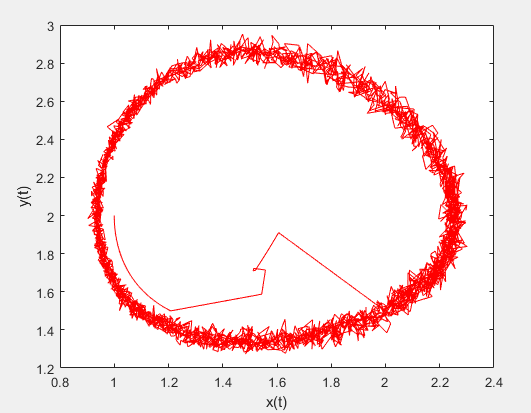}}\\
		\caption*{}
		\label{fig:fig12b}
		\subfloat[for $\delta=0.87$]{\includegraphics[width=0.45\textwidth]{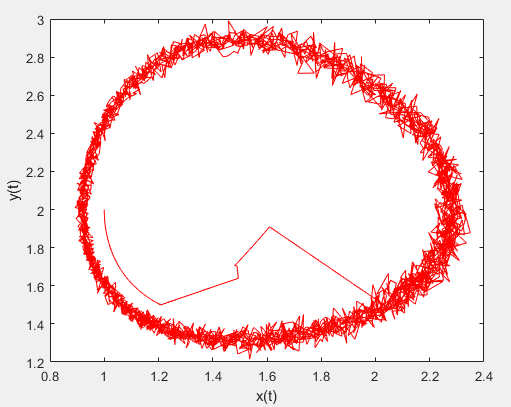}}\hspace{0.5cm}
		\label{fig:fig12c}
		\subfloat[for $\delta=0.98$]{\includegraphics[width=0.47\textwidth]{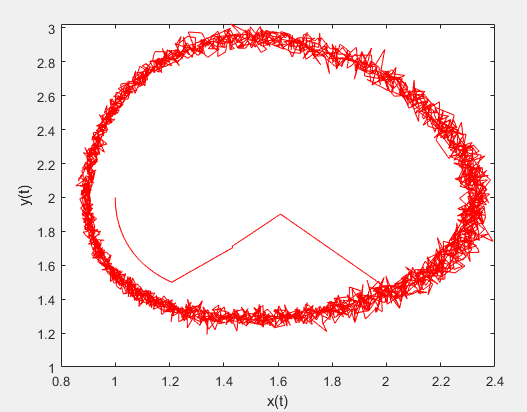}}\\
		%		\caption*{B-picture (a) and (b)}
		\caption{Phase Portrait between $x(t)$ and $y(t)$} 
		\label{fig:fig12d}
	\end{figure}

	By using the exponenetial decay, it presents the numerical solutions to the piecewise Lotka Volterra predator-prey model.The starting circumstances were regarded as
	$t(0)=0$,  $h=0.01$, $x(0)=1$, $y(0)=2$, $\lambda_{1}=2$, $\lambda_{2}=1$, $\lambda_{3}=1.5$, $\lambda_{4}=1$, $\sigma_{1}=0.01$ and $\sigma_{2}=0.02$. 
	See figures \ref{fig:fig12d}(a), \ref{fig:fig12d}(b), \ref{fig:fig12d}(c) and \ref{fig:fig12d}(d).

	\section{Conclusion} 
	One of the key connections that controls ecosystem yields is that between prey and predator.   Each species' Lotka Volterra predator model is studied and investigated in this study. The piecewise system of the connections between predator and prey is significant to the investigation. 
	The theory behind the piecewise system of the connections in the model has advanced considerably. Additionally, the model's stability and uniqueness were also the topics of the inquiry. 
	The examination displays the link between the predator and prey. The work presents the idea of piecewise derivatives for different fractional derivatives in the modeling of the interaction between prey and predator. 
	Computer simulations were used to explain and interpret the numerical results. The computer results show the real-issues behaviors of the predator-prey relationship modeling. For natural populations, none of the aforementioned presumptions are likely to be true. However, the Lotka–Volterra model reveals two crucial characteristics of populations of predators and prey, and these characteristics frequently hold true for model variations where these presumptions are relaxed.

\end{document}